\documentclass[12pt]{amsart}
\usepackage{amsmath}
\usepackage{amsfonts}
\usepackage{amssymb}
\usepackage{amscd}
\usepackage{graphicx}
\usepackage[abbrev,alphabetic]{amsrefs}
\RequirePackage[dvipsnames,usenames]{color}
\usepackage{soul,xcolor}
\setstcolor{red}
\usepackage{stmaryrd}
\usepackage{mathtools}
\usepackage{booktabs}
\usepackage{multirow}
\usepackage{stmaryrd} 
\newtagform{tiny}{\tiny(}{)}

\usepackage{mathtools}
\usepackage{hyperref}
\usepackage[margin=1.25in]{geometry}

\usepackage{amsthm}
\usepackage{comment}
\usepackage[all,cmtip]{xy}
\usepackage{tikz-cd}
\usetikzlibrary{cd}
\usepackage{mathtools}
\usepackage[all]{xy}

\makeatletter
\def\@tocline#1#2#3#4#5#6#7{\relax
  \ifnum #1>\c@tocdepth 
  \else
    \par \addpenalty\@secpenalty\addvspace{#2}%
    \begingroup \hyphenpenalty\@M
    \@ifempty{#4}{%
      \@tempdima\csname r@tocindent\number#1\endcsname\relax
    }{%
      \@tempdima#4\relax
    }%
    \parindent\z@ \leftskip#3\relax \advance\leftskip\@tempdima\relax
    \rightskip\@pnumwidth plus4em \parfillskip-\@pnumwidth
    #5\leavevmode\hskip-\@tempdima
      \ifcase #1
       \or\or \hskip 1em \or \hskip 2em \else \hskip 3em \fi%
      #6\nobreak\relax
    \hfill\hbox to\@pnumwidth{\@tocpagenum{#7}}\par
    \nobreak
    \endgroup
  \fi}
\makeatother

\newsavebox{\pullback}
\sbox\pullback{%
\begin{tikzpicture}%
\draw (0,0) -- (1ex,0ex);%
\draw (1ex,0ex) -- (1ex,1ex);%
\end{tikzpicture}}

\newsavebox{\pullbackdl}
\sbox\pullbackdl{%
\begin{tikzpicture}%
\draw (-1ex,0ex) -- (0ex,0ex);%
\draw (0ex,-1ex) -- (0ex,0ex);%
\end{tikzpicture}}

\newsavebox{\pushoutdr}
\sbox\pushoutdr{%
\begin{tikzpicture}%
\draw (-1ex,-1ex) -- (-1ex,0ex);%
\draw (-1ex,0ex) -- (0ex,0ex);%
\end{tikzpicture}}

\newcommand{\rup}[1]{\lceil #1 \rceil}
\newcommand{\rdown}[1]{\lfloor #1 \rfloor}

\renewcommand{\P}{\mathbb{P}}
\newcommand{\Z}{\mathbb{Z}}
\newcommand{\Q}{\mathbb{Q}}

\newcommand{\R}{\mathbb{R}}
\newcommand{\F}{\mathbb{F}}

\newcommand{\cO}{\mathcal{O}}

\newcommand{\MO}{\mathcal{O}}
\newcommand{\sO}{\mathcal{O}}

\newcommand{\m}{\mathfrak{m}}

\DeclareMathOperator{\pr}{pr}

\DeclareMathOperator{\Supp}{Supp}
\DeclareMathOperator{\Spec}{Spec}

\DeclareMathOperator{\Hom}{Hom}

\DeclareMathOperator{\Ex}{Ex}

\DeclareMathOperator{\Coker}{Coker}
\DeclareMathOperator{\Bl}{Bl}

\theoremstyle{plain}
\newtheorem{theorem}{Theorem}[section]
\newtheorem{thm}[theorem]{Theorem}

\newtheorem{prop}[theorem]{Proposition}

\newtheorem{lem}[theorem]{Lemma}

\newtheorem{cor}[theorem]{Corollary}

\newtheorem*{claim*}{Claim}

\theoremstyle{definition}

\newtheorem{dfn}[theorem]{Definition}

\newtheorem{Convention}[theorem]{Convention}
\newtheorem{example}[theorem]{Example}

\newtheorem*{setup*}{Setup}
\theoremstyle{plain}
\newtheorem{theo}{Theorem}

\theoremstyle{remark}
\newtheorem{rem}[theorem]{Remark}



\numberwithin{equation}{theorem}

\title[Global $F$-regularity for weak del Pezzo surfaces]{Global $F$-regularity for weak del Pezzo surfaces}

\author{Tatsuro Kawakami and Hiromu Tanaka}
\address{Department of Mathematics, Graduate School of Science, Kyoto University, Kyoto 606-8502, Japan} 
\email{tatsurokawakami0@gmail.com}

\address{Graduate School of Mathematical Sciences, 
The University of Tokyo, 
3-8-1 Komaba, Meguro-ku, Tokyo 153-8914, JAPAN} 
\email{tanaka@ms.u-tokyo.ac.jp}

\begin{document}

\begin{abstract}
Let $k$ be an algebraically closed field of characteristic $p>0$.
Let $X$ be a normal projective surface over $k$ with canonical singularities whose anti-canonical divisor is nef and big.
We prove that $X$ is globally $F$-regular except for the following cases:
(1) $K_X^2=4$ and $p=2$, 
(2) $K_X^2=3$ and $p \in \{2, 3\}$,  
(3) $K_X^2=2$ and $p \in \{2, 3\}$,  
(4) $K_X^2=1$ and $p \in \{2, 3, 5\}$. 
For each degree $K_X^2$, the assumption of $p$ is optimal.
\end{abstract}

\subjclass[2020]{14J45, 13A35}   
\keywords{del Pezzo surfaces, globally $F$-regular}
\maketitle

\setcounter{tocdepth}{2}

\tableofcontents

\section{Introduction}

We work over an algebraically closed field of characteristic $p>0$.
Fano varieties play a significant role in the classification of algebraic varieties.
In positive characteristic, properties defined by the Frobenius morphism such as (global) $F$-splitting or global $F$-regularity are useful.
Therefore, it is natural to ask when Fano varieties are $F$-split or globally $F$-regular.
For smooth del Pezzo surfaces, Hara \cite{Hara98} proved the following result:

\begin{thm}[\textup{\cite{Hara98}*{Example 5.5}}]
    Let $X$ be a smooth del Pezzo surface over an algebraically closed field of characteristic $p>0$.
    Then $X$ is globally $F$-regular except for the following cases:
    \begin{enumerate}
        \item $K_X^2=3$ and $p=2$.
        \item $K_X^2=2$ and $p \in \{2, 3\}$. 
        \item $K_X^2=1$ and $p \in \{2, 3, 5\}$. 
    \end{enumerate}
\end{thm}

Our aim is to generalize Hara's result to the case when $-K_X$ is nef and big, or equivalently, $X$ is canonical.
Here, we say that a variety is \textit{canonical} if it has only canonical singularities.
In fact, the following theorem holds.

\begin{theo}\label{Intro:main}
    Let $k$ be an algebraically closed field of characteristic $p>0$.
    Let $X$ be a canonical projective surface over $k$ whose anti-canonical divisor is nef and big. 
    Then $X$ is globally $F$-regular except for the following cases:
    \begin{enumerate}
    \item $K_X^2=4$ and $p=2$.
    \item $K_X^2=3$ and  $p \in \{2, 3\}$. 
    \item $K_X^2=2$ and  $p \in \{2, 3\}$. 
    \item $K_X^2=1$ and  $p \in \{2, 3, 5\}$. 
    \end{enumerate}
\end{theo}

\noindent Theorem \ref{Intro:main} has an important role for investigating of global $F$-regularity of smooth del Pezzo varieties and smooth Fano threefolds (see \cites{Kawakami-Tanaka(dPvar),Kawakami-Tanaka(Fano)} for details).

\begin{rem}
    For each degree $K_X^2$, the assumption on $p$ is optimal.
    In fact, there exists a canonical del Pezzo surface that is not strongly $F$-regular for each listed above (see \cite{KN2}*{Table 1} and \cite{Kawakami-Takamatsu}*{Table 1}).
    Moreover, the assumption on $p$ is still optimal even if we assume that $X$ is smooth, since taking the minimal resolution does not change the degree and globally $F$-regularity (see Section \ref{subsubsection:weak dP} and Corollary \ref{c-F-split-blowup}).
\end{rem}

\begin{rem}
    For each prime number $p$, there exists a klt del Pezzo surface $X$ 
    (i.e.,  a normal projective surface such that $(X, 0)$ is klt and $-K_X$ is ample) that is not $F$-split (\cite{CTW2}*{Theorem 1.1}). 
\end{rem}

We now focus on the proof of Theorem \ref{Intro:main}.
We first investigate when $X$ as in Theorem \ref{Intro:main} is $F$-split. 
The proof is divided into two cases: the case where $K_X^2\geq 5$ and the case where $K_X^2\leq 4$.

First, we consider the case where $K_X^2\geq 5$. 
We may assume that $X$ 
is obtained by taking a blowup $f\colon X \to \P^2$ along some points. 
Recall that 
if there exists an effective divisor $\Delta_{\P^2}$ on $\P^2$ 
such that the divisor $\Delta$ on $X$ defined by 
$K_X+\Delta = f^*(K_{\P^2}+\Delta_{\P^2})$ is effective, 
then the following holds (Proposition \ref{p-F-split-blowup}): 
\[
(\P^2, \Delta_{\P^2}) \text{ is $F$-split} \Leftrightarrow 
(X, \Delta) \text{ is $F$-split} \Rightarrow X \text{ is $F$-split}. 
\]
Such a divisor $\Delta_{\P^2}$ can be found by utilizing 
an inversion of adjunction for $F$-splitting (Proposition \ref{p-F-split-IOA}). 
However, since we can only assume that the blowup points are in \textit{almost general position}, the situation is more complicated 
than the smooth del Pezzo cases, which are obtained by blowing up points in \textit{general position}. 
For more details, see Proposition \ref{prop:degree 5}.

Next, we consider the case where $K_X^2\leq 4$. In the case of smooth del Pezzo surfaces, Hara \cite{Hara98} investigated $F$-splitting using the following two steps:

\begin{enumerate}
     \item[(i)] The reduction of $F$-splitting to the vanishing $H^1(X, \Omega_X^{1}(pK_X))$ via the Cartier operator.
    \item[(ii)]  Embedding $X$ as a hypersurface or a complete intersection of a weighted projective space and proving the vanishing $H^1(X, \Omega_X^{1}(pK_X))$.
\end{enumerate}

For a smooth weak del Pezzo surface $X$, we can also reduce the $F$-splitting of $X$ to the vanishing $H^1(X, \Omega_X^{1}(pK_X))$.
However, Step (ii), proving $H^1(X, \Omega_X^{1}(pK_X))=0$, is not easy because $X$ is not embedded in a weighted projective space via $|-mK_X|$ for $m\in\Z_{>0}$. 
Therefore, by replacing $X$ with its anti-canonical model, we embed $X$ into a weighted projective space via $|-mK_X|$. 
However, in this case, Step (i), the reduction of the $F$-splitting of $X$ to the vanishing of $H^1(X, \Omega_X^{1}(pK_X))$, is not straightforward since Cartier operator is defined on smooth schemes.
To address this issue, we use the reflexive Cartier operator introduced in \cite{Kaw4}.
Indeed, we utilize the fact that the reflexive Cartier operator behaves well on $F$-pure klt surfaces (Lemma \ref{lem:exact sequence}).

Combining the above results, we obtain the following theorem.

\begin{theo}\label{Intro:F-split}
    Let $k$ be an algebraically closed field of characteristic $p>0$.
    Let $X$ be a canonical $F$-pure projective surface over $k$ 
     whose anti-canonical divisor is ample. 
    Then $X$ is $F$-split except for the following cases:
    \begin{enumerate}
        \item $K_X^2=3$ and $p=2$.
        \item $K_X^2=2$ and  $p \in \{2, 3\}$. 
        \item $K_X^2=1$ and  $p \in \{2, 3, 5\}$. 
    \end{enumerate}
\end{theo}
\begin{rem}\,
    \begin{enumerate}
        \item When $K_X^2=4$ and $p=2$, there exists a a canonical del Pezzo surface with $D^0_5$-singularity, which is not $F$-pure (\cite{KN2}*{Proposition 3.16}).
        \item When $K_X^2=3$ and $p=3$, there exists a canonical del Pezzo surface with $E^0_6$-singularity, which is not $F$-pure (\cite{KN2}*{Proposition 3.22}). 
    \end{enumerate}
\end{rem}

Finally, we now overview 
how to deduce Theorem \ref{Intro:main} from Theorem \ref{Intro:F-split}. Let 
$X$ be a canonical weak del Pezzo surface. 
Replacing $X$ with its anti-canonical model, we can assume that $-K_X$ is ample.
For each degree $K_X^2$, we will find an optimal bound on $p$ that ensures all the singularities of $X$ are strongly $F$-regular (see Lemma \ref{lem:F-reg-sing}).
We then conclude by the equivalence between $F$-splitting and global $F$-regularity for strongly $F$-regular $\Q$-Gorenstein Fano varieties  (cf.~\cite{Kawakami-Totaro}*{Lemma 6.4}).

\medskip
\noindent {\bf Acknowledgements.}
The authors express their gratitude to Burt Totaro for valuable comments.
Kawakami was supported by JSPS KAKENHI Grant number JP22J00272.
Tanaka was supported by JSPS KAKENHI Grant number JP22H01112 and JP23K03028. 

\section{Preliminaries}

\subsection{Notation and terminology}
\label{ss:notation}

In this subsection, we summarize notation and basic definitions used in this article. 
\begin{enumerate}
\item Throughout the paper, $p$ denotes a prime number and 
we work over an algebraically closed field $k$ of characteristic $p>0$. 
We set $\F_p\coloneqq\Z/p\Z$. 
We denote by $F \colon X \to X$ the absolute Frobenius morphism on an $\F_p$-scheme $X$.
    \item We say that $X$ is a {\em variety} (over $k$) if 
    $X$ is an integral scheme 
    that is separated and of finite type over $k$. 
    We say that $X$ is a {\em curve} (resp.~{\em surface}) 
    if $X$ is a variety of dimension one (resp.~two). 
\item For a variety $X$, 
we define the {\em function field} $K(X)$ of $X$ 
as the stalk $\MO_{X, \xi}$ at the generic point $\xi$ of $X$. 
\item We say that a $\Q$-divisor $D$ on a normal variety $X$ is \emph{simple normal crossing} if for every point $x \in \Supp D$, the local ring $\cO_{X,x}$ is regular and there exists a regular system of parameters $x_1,\ldots, x_d$ of the maximal ideal $\m$ of $\cO_{X,x}$ and $1 \leq r \leq d$ such that $\Supp (D|_{\Spec \MO_{X, x}}) = \Spec (\cO_{X,x}/(x_1 \cdots x_r))$. 
\item Given a variety $X$, a projective birational morphism $\pi \colon Y \to X$ is called a \emph{log resolution of $X$} if $Y$ is a smooth variety and $\mathrm{Exc}(f)$ is a simple normal crossing divisor.
\item Given a variety $X$ and a closed subscheme $Z$, 
we denote by $\Bl_Z X$ the blowup of $X$ along $Z$. 
\item Given a normal variety $X$ and a $\Z$-divisor $D$ on $X$, we define a reflexive sheaf $\Omega_X^{[i]}(D)$ by $j_{*}(\Omega_U^{i}\otimes \sO_U(D))$, where $j\colon U\hookrightarrow X$ is the open immersion from the smooth locus 
$U$ of $X$.
\end{enumerate}

\subsubsection{Singularities in minimal model program}
For the definitions of singularities in minimal model program (e.g., canonical and klt), we refer to \cite{KM98}*{Section 2.3}. 
Take a normal surface $X$. 
Let $f\colon Y \to X$ be the minimal resolution. We only need the following characterizations in this paper. 
\begin{enumerate}
\item 
$X$ is canonical if and only if    $K_X$ is $\Q$-Cartier and $K_Y \sim f^*K_X$. 
\item 
$X$ is klt if and only if  $K_X$ is $\Q$-Cartier, $f$ is a log resolution of $X$, and 
all the coefficients of the $\Q$-divisor $\Gamma$ defined by $K_Y+\Gamma \sim f^*K_X$ are $<1$.  
\end{enumerate}
By definition, we have 
\[
\text{canonical}\Rightarrow 
\text{klt}. 
\]
Moreover, the following implications hold for the surface case: 
\begin{itemize}
\item canonical $\Rightarrow$ Gorenstein. 
\item klt $\Rightarrow$ $\Q$-factorial. 
\end{itemize}

\subsubsection{Weak del Pezzo surfaces}\label{subsubsection:weak dP}

Given a normal projective Gorenstein surface $X$, 
we say that $X$ is {\em del Pezzo} (resp.~{\em weak del Pezzo}) 
if $-K_X$ is ample (resp.~nef and big). 
In what follows, we summarize some properties on (weak) del Pezzo surfaces for later usage.

Let $Z$ be a canonical weak del Pezzo surface. 
The {\em anti-canonical model} $Y$ of $Z$ is defined 
as the Stein factorisation of the morphism $\varphi_{|-mK_Z|} \colon Z \to \P^{h^0(Z, -mK_Z)-1}$ induced by the complete linear system $|-mK_Z|$, 
where $m$ is a positive integer such that $|-mK_X|$ is base point free (whose existence is guaranteed by \cite{Tan15}*{Theorem 0.4}). 
Then $Y$ is canonical, because  have $K_Z \sim h^*K_Y$ for the induced morphism $h \colon Z \to Y$. 
Moreover, $h$ is obtained by contracting all the $(-2)$-curves on $Y$. 
In particular, the minimal resolution $X$ of $Z$ coincides with the minimal resolution of $Y$: 
\[
f \colon X \xrightarrow{g} Z \xrightarrow{h} Y. 
\]
Moreover, $X$ is a smooth weak del Pezzo surface. 

There is a natural one-to-one correspondence between
\begin{itemize}
\item smooth weak del Pezzo surfaces, and 
\item canonical del Pezzo surfaces. 
\end{itemize}
Indeed, if $Y$ is a canonical del Pezzo surface, then its minimal resolution $X$ is a smooth weak del Pezzo surface. 
Conversely, given a smooth weak del Pezzo surface $X$, 
its anti-canonical model $Y$ is a canonical del Pezzo surface.

\subsection{\texorpdfstring{$F$}--splitting and global \texorpdfstring{$F$}--regularity}

In this subsection, we gather basic facts on $F$-splitting and global $F$-regularity. 

\begin{dfn}
Let $X$ be a normal variety and let $\Delta$ be an effective $\Q$-divisor on $X$. 
\begin{enumerate}
\item We say that $(X, \Delta)$ is {\em $F$-split} if 
\[
\MO_X \xrightarrow{F^e} F_*^e\MO_X \hookrightarrow F_*^e\MO_X( \rdown{(p^e-1)\Delta})
\]
splits as an $\MO_X$-module homomorphism for every $e \in \Z_{>0}$. 
\item We say that $(X, \Delta)$ is {\em sharply $F$-split} if 
\[
\MO_X \xrightarrow{F^e} F_*^e\MO_X \hookrightarrow F_*^e\MO_X( \rup{(p^e-1)\Delta})
\]
splits as an $\MO_X$-module homomorphism for some $e \in \Z_{>0}$. 
\item 
We say that $(X, \Delta)$ is {\em globally $F$-regular} 
if, given an effective $\Z$-divisor $E$, there exists $e \in \Z_{>0}$ such that 
\[
\MO_X \xrightarrow{F^e} F_*^e\MO_X \hookrightarrow 
F_*^e\MO_X( \rup{(p^e-1)\Delta} +E)
\]
splits as an $\MO_X$-module homomorphism. 
\end{enumerate}
We say that $X$ is {\em $F$-split} (resp. {\em globally $F$-regular}) 
if so is $(X, 0)$. 
\end{dfn}

\begin{rem}
We have the following implications: 
\[
\text{globally $F$-regular}\Longrightarrow
\text{sharply $F$-split}\Longrightarrow
\text{$F$-split} 
 \]
 where the former implication is easy and the latter one holds by 
 the same argument as in  \cite{Sch08g}*{Proposition 3.3}. 
 Moreover, 
 if the condition ($\star$) holds, then 
 $(X, \Delta)$ is sharply F-split if and only if $(X, \Delta)$ is $F$-split. 
\begin{enumerate}
    \item[($\star$)] $(p^e-1)\Delta$ is a $\Z$-divisor for some $e \in \Z_{>0}$. 
 \end{enumerate}
In particular, $X$ is $F$-split if and only if $F\colon\MO_X \to F_*\MO_X$ splits as an $\MO_X$-module homomorphism. 
 In this paper, we only treat the case when  ($\star$) holds, and hence being $F$-split is equivalent to being sharply $F$-split. 
 For more foundational properties, we refer to \cite{SS10}. 
\end{rem}

We shall also use the local versions of $F$-splitting and global $F$-regularity. 

\begin{dfn}
Given a normal variety $X$, 
we say that $X$ is {\em $F$-pure} (resp.~{\em strongly $F$-regular}) 
if there exists an open cover $X = \bigcup_{i \in I} X_i$ such that 
$X_i$ is {$F$-split} (resp.~{globally $F$-regular}) 
for every $i \in I$. 
\end{dfn}

In what follows, we summarize some $F$-splitting criteria, which are well known to experts.

\begin{prop}\label{p-F-split-blowup}
Let $f: X \to Y$ be a birational morphism of normal projective varieties. 
Take an effective $\Q$-divisor $\Delta_Y$ on $Y$ such that 
$(p^e-1)(K_Y+\Delta_Y)$ is Cartier for some $e \in \Z_{>0}$. 
Assume that the $\Q$-divisor $\Delta$ defined by $K_X+\Delta = f^*(K_Y+\Delta_Y)$ is effective. 
Then $(X, \Delta)$ is  $F$-split (resp.\ globally $F$-regular) if and only if $(Y, \Delta_Y)$ is $F$-split (resp.\ globally $F$-regular).  
\end{prop}

\begin{proof}
If $(X, \Delta)$ is $F$-split, then so is $(Y, \Delta_Y)$, 
which can be checked by taking the push-forward. 
As for the opposite implication, the same argument as in 
\cite{HX15}*{the first paragraph of the proof of Proposition 2.11} works. 
\end{proof}

\begin{cor}\label{c-F-split-blowup}
Let $Y$ be a canonical projective surface and 
let $f: X \to Y$ be its minimal resolution. 
Then $X$ is $F$-split (resp.\ globally $F$-regular)  if and only if $Y$ is $F$-split (resp.\ globally $F$-regular).  
\end{cor}

\begin{proof}
The assertion immediately follows from Proposition \ref{p-F-split-blowup} 
by using $K_X =f^*K_Y$. 
\end{proof}

\begin{prop}\label{p-F-split-IOA}
Let $X$ be a normal projective Gorenstein variety. 
Take a normal prime Cartier divisor $S$ 
and an effective $\Q$-Cartier $\Q$-divisor $B$ on $X$ such that $S \not\subset \Supp\,B$. 
Assume that 
\begin{enumerate}
\item $(S, B|_S)$ is $F$-split, and 
\item there is a positive integer $e \in \Z_{>0}$  such that 
$(p^e-1)(K_X+S+B)$ is Cartier and 
\[
H^1(X, \MO_X(-S-(p^e-1)(K_X+S+B)))=0. 
\]
\end{enumerate}
Then $(X, S+B)$ is $F$-split. 
\end{prop}

\begin{proof}
The same argument as in \cite{CTW17}*{Lemma 2.7} works. 
\end{proof}

\begin{example}\label{e-P2-P1P1}
We now summarize some easy cases for later usage, 
although all of them are well known to experts, 
\begin{enumerate}
\item If $P, Q \in \P^1$ are distinct points, then $(\P^1, P+Q)$ is $F$-split (Proposition \ref{p-F-split-IOA}). 
\item Let $L, L', L''$ be three lines on $\P^2$ such that $L+L'+L''$ is simple normal crossing. 
Then $(\P^2, L+L'+L'')$ is $F$-split by (1) and Proposition \ref{p-F-split-IOA}. 
\item For each $i \in \{1, 2\}$, let $F_i$ and $F'_i$ be distinct fibers 
of the $i$-th projection $\pr_i \colon \P^1 \times \P^1 \to \P^1$. 
Then $(\P^1 \times \P^1, F_1 +F'_1 + F_2+F'_2)$ is $F$-split by (1) and Proposition \ref{p-F-split-IOA}. 
\end{enumerate}
\end{example}

\subsection{Reflexive Cartier operator}

Throughout this subsection, we use the following convention unless stated otherwise.

\begin{Convention}
Let $X$ be a normal variety and $D$ a $\Z$-divisor on $X$.
Let $U$ be the smooth locus of $X$ and $j\colon U\hookrightarrow X$ the inclusion. 
By abuse of notation, $D|_U$ is denoted by $D$. 
\end{Convention}

The Frobenius pushforward of the de Rham complex
\[
F_{*}\Omega^{\bullet}_U\colon  F_{*}\sO_U \xrightarrow{F_{*}d} F_{*}\Omega^1_U \xrightarrow{F_{*}d} \cdots
\]
is a complex of $\sO_U$-modules.
Tensoring with $\sO_U(D)$, we obtain a complex 
\[
F_{*}\Omega^{\bullet}_U\colon  F_{*}\sO_U(pD) \xrightarrow{F_{*}d\otimes \sO_U(D)} F_{*}\Omega^1_U(pD) \xrightarrow{F_{*}d\otimes \sO_U(D)} \cdots
\]

We define coherent $\sO_U$-modules $B^{i}_U(pD)$ and $Z^{i}_U(pD)$ by 
\begin{align*}
    &B^{i}_U(pD)\coloneqq \mathrm{Im}(F_{*}d\otimes\sO_U(D) \colon F_{*}\Omega^{i-1}_U(pD) \to F_{*}\Omega^{i}_U(pD)),\\
    &Z^{i}_U(pD)\coloneqq \mathrm{Ker}(F_{*}d\otimes\sO_U(D) \colon F_{*}\Omega^{i}_U(pD) \to F_{*}\Omega^{i+1}_U(pD)),
\end{align*}
for all $i\geq 0$.
Then $B^{i}_U(pD)$ and $Z^{i}_U(pD)$ are locally free (\cite{Kaw4}*{Lemma 3.2}).
{When $D=0$, we simply denote $B^{i}_U(pD)$ and $Z^{i}_U(pD)$ by $B^{i}_U$ and $Z^{i}_U$, respectively.
Then $B^{i}_U(pD)=B^{i}_U\otimes \sO_U(D)$ and $Z^{i}_U(pD)=Z^{i}_U\otimes \sO_U(D)$ holds (\cite{Kaw4}*{Remark 3.3}). In particular, we note that $B^{i}_U(pD)$ and $Z^{i}_U(pD)$ do not mean $B^{i}_U\otimes \sO_U(pD)$ and $Z^{i}_U\otimes \sO_U(pD)$.}

By \cite{Kaw4}*{Lemma 3.2}, there exists an exact sequence
\begin{equation}
    0 \to B^{i}_U(pD)\to Z^{i}_U(pD)\xrightarrow{C^{i}_{U}(D)} \Omega^{i}_U(D)\to 0,\label{log smooth 2}
\end{equation}
and the map $C^{i}_{U}(D)$ coincides with $C^{i}_{U}\otimes \sO_U(D)$, where  $C^{i}_{U}$ is the usual Cartier operator.

\begin{dfn}\label{def:reflexive Carter operators}
We define reflexive $\sO_X$-modules 
$B^{[i]}_X(pD)$ and $Z^{[i]}_X(pD)$ by 
\begin{align*}
    &B^{[i]}_X(pD)\coloneqq j_{*}B^{i}_U(pD)\,\,\text{and}\\
    &Z^{[i]}_X(pD)\coloneqq j_{*}Z^{i}_U(pD)
\end{align*}
for all $i\geq 0$.
The \textit{$i$-th reflexive Cartier operator}
\[
C^{[i]}_{X}(D)\colon Z^{[i]}_X(pD)\to \Omega^{[i]}_X(D)
\]
\textit{associated to $D$} is defined as $j_{*}C^{i}_{U}(D)$ for all $i\geq 0$.
\end{dfn}

\begin{lem}\label{lem:Cartier operators}
There exist the following exact sequences:
\begin{equation}
0 \to Z^{[i]}_X(pD) \to F_{*}\Omega^{[i]}_X(pD) \xrightarrow{d'} B^{[i+1]}_X(pD),\label{non-log smooth 1}
\end{equation}
\begin{equation}
0 \to B^{[i]}_X(pD)\to Z^{[i]}_X(pD)\xrightarrow{C^{[i]}_{X}(D)} \Omega^{[i]}_X(D),\label{non-log smooth 2}
\end{equation}
for all $i\geq 0$.
Moreover, 
$d'|_U \colon F_{*}\Omega^{[i]}_X(pD)|_U \to  B^{[i+1]}_X(pD)|_U$ and 
$C^{[i]}_{X}(D)|_U \colon Z^{[i]}_X(pD)|_U  \to \Omega^{[i]}_X(D)|_U$ are surjective,
and the homomorphism $C^{[i]}_{X}(D)|_{U}$ coincides with $C^{i}_{U}\otimes\sO_U(D)$. 
\end{lem}
\begin{proof}
Taking $B=0$ in \cite{Kaw4}*{Lemma 3.5}, we obtain the assertion. 
\end{proof}

\begin{rem}\label{rem:B^1}
Taking $i=0$ in \eqref{non-log smooth 1}, we obtain an exact sequence
\[
0 \to \sO_X(D) \to F_{*}\sO_X(pD) \to B_X^{[1]}(pD),
\]
and the first map is induced by the Frobenius homomorphism.
In particular, \[B_X^{[1]}(pD)=j_{*}\Coker({F\colon \sO_U(D) \to F_{*}\sO_U(pD)})\]
holds.
\end{rem}

\section{Proofs of main theorems}

\subsection{Criterion of the \texorpdfstring{$F$}--splitting of klt surfaces}

In this subsection, we provide a criterion for the $F$-splitting of klt surfaces (Proposition \ref{prop:GFR criterion}).

\begin{lem}\label{lem:exact sequence}
    Let $X$ be an $F$-pure klt surface and $D$ a $\Z$-divisor.
    Then the  sequence 
    \begin{align}
0 \to B^{[i]}_X(pD) \to  Z^{[i]}_X(pD) \xrightarrow{C_X^{[i]}(D)} \Omega^{[i]}_X(D) \to 0 \label{exact:Bsingular}
\end{align}
 is exact.
\end{lem}

\begin{proof}
It is enough to show that $C_X^{[i]}(D)$ is surjective, as the other parts has been settled in Lemma \ref{lem:Cartier operators}.
    Since the assertion is local on $X$, we may assume that $X$ is affine and has a unique singular point $Q$. 
    If $p\neq 5$ or the singularity $Q$ is not RDP of type $E_8^1$, then $X$ is $F$-liftable by \cite{Kawakami-Takamatsu}*{Theorem A}. Then the surjectivity of $C^{[i]}(D)$ follows from \cite{Kaw4}*{Lemma 3.8}.

    Suppose that $p=5$ and the singularity  $Q$ is of type $E_8^1$.
    Then we may assume that $D=0$ by \cite{Lipman69}*{Section 24} (see also \cite{LMM2}*{Table 2}). Then the desired surjectivity follows from \cite{Kaw4}*{Proposition 4.4} and \cite{Kawakami-Takamatsu}*{Theorem B}.
\end{proof}

\begin{prop}\label{prop:GFR criterion}
    Let $X$ be an $F$-pure klt projective surface.
    Suppose that the following conditions hold:
    \begin{enumerate}
        \item $H^0(X, \Omega_X^{[1]}(K_X))=0$.
        \item $H^1(X, \Omega_X^{[1]}(pK_X))=0$.
        \item $H^0(X, \sO_X((p+1)K_X))=0$.
    \end{enumerate}
    Then $X$ is globally $F$-split.
\end{prop}
\begin{proof}
Recall that $X$ is $F$-split if and only if 
the evaluation map
    \[
    \mathrm{Hom}_{\sO_X}(F_{*}\sO_X, \sO_X) \xrightarrow{F^*} \Hom_{\MO_X}(\MO_X, \MO_X) (\cong  H^0(X, \sO_X)) 
    \]
    is surjective. 
    By Serre duality, they are also equivalent to the injectivity of 
    \begin{equation}
        H^2(X, \sO_X(K_X))\to H^2(X, \sO_X(pK_X)).\label{eq:inj}
    \end{equation}

    Let $U$ be the smooth locus of $X$.
    Since $X$ is $F$-pure, the exact sequence 
    \[
    0\to \sO_U \to F_{*}\sO_U \to B_U^{1}\to 0
    \]
    splits locally.
    Tensoring with $\sO_U(K_X)$ and taking the pushforward by the inclusion $U\hookrightarrow X$, we have the following locally split exact sequence
    \[
    0\to \sO_X(K_X) \to F_{*}\sO_X(pK_X) \to B^{[1]}_X(pK_X)\to 0
    \]
    by Remark \ref{rem:B^1}.

    Thus, for the injectivity of \eqref{eq:inj}, it suffices to show that $H^1(X, B^{[1]}_X(pK_X))=0$.
    By 
    Lemma \ref{lem:exact sequence} and the condition (1), 
    it is enough to prove that 
    $H^1(X, Z_X^{[1]}(pK_X))=0$.
    By \eqref{non-log smooth 1} and the condition (2),
    it suffices to show that 
    \[
    H^0(X, B_X^{[2]}(pK_X))=0.
    \]
    Since $B_X^{[2]}(pK_X)$ is a subsheaf of $F_{*}\Omega_X^{[2]}(pK_X)$, 
    we have
    \[
    H^0(X, B_X^{[2]}(pK_X))\hookrightarrow H^0(X, \Omega_X^{[2]}(pK_X))=H^0(X, \sO_X((p+1)K_X)).
    \]
    By the condition (3), we conclude.
\end{proof}

The condition (3) of Proposition \ref{prop:GFR criterion}  is satisfied 
if $-K_X$ is big.
In what follows, we see when the condition (1) of Proposition \ref{prop:GFR criterion} is satisfied.

\begin{dfn}[Log liftability]
Let $X$ be a normal projective surface.
We say that $X$ is \textit{log liftable} if there exists a log resolution $f\colon Y \to X$ of $X$ such that $(Y,\mathrm{Exc}(f))$ lifts to the ring $W(k)$ of Witt vectors.
For the definition of liftability of a log smooth pair, we refer to \cite{Kaw3}*{Definition 2.6}.
\end{dfn}

\begin{lem}\label{lem:log lift vs H^0}
    Let $X$ be a normal projective $F$-pure surface 
    such that $-K_X$ is a nef and big $\Q$-Cartier divisor.
    Then $X$ is log liftable if and only if $H^0(X, \Omega_X^{[1]}(K_X))=0$.
\end{lem}
\begin{proof}
    Since $H^2(X, \sO_X)\cong H^0(X, \sO_X(K_X))=0$, the `if' direction is \cite{Kaw5}*{Theorem 2.8}.
    We prove the `only if' direction. 
    Let $f\colon Y\to X$ be a log resolution such that $(Y,E\coloneqq \Ex(f))$ lifts to $W(k)$.
    Since $f_{*}(\Omega^1_Y(\log\,E)(f^{*}K_X))=\Omega_X^{[1]}(K_X)$ by \cite{Kawakami-Takamatsu}*{Theorem B}, we have $H^0(X, \Omega_X^{[1]}(K_X))=H^0(Y,\Omega_Y^{1}(\log\,E)\otimes \sO_Y(f^{*}K_X))$.
    Then the vanishing follows from \cite{Kaw3}*{Theorem 2.11}.
\end{proof}

\subsection{Global \texorpdfstring{$F$}--splitting: Proof of Theorem \ref{Intro:F-split}}

In the following proposition, we investigate $F$-splitting of $F$-pure canonical del Pezzo surfaces. 
For the proof, we confirm when the condition (2) of Proposition \ref{prop:GFR criterion} is satisfied.

\begin{prop}\label{prop:degree 1--4}
Let $X$ be an $F$-pure canonical del Pezzo surface. 
Suppose that one of the following holds.
\begin{enumerate}
    \item $K_X^2=1$ and $p>5$.
    \item $K_X^2=2$ and $p>3$.
    \item $K_X^2=3$ and $p>2$.
    \item $K_X^2=4$.
\end{enumerate}
Then $X$ is $F$-split.
\end{prop}
\begin{proof}
In each case, $X$ is log liftable by \cite{KN}*{Theorem 1.7 (1)}, and thus the condition (1) of Proposition \ref{prop:GFR criterion} is satisfied.
Thus it suffices to confirm the condition (2) of Proposition \ref{prop:GFR criterion}, i.e., $H^1(X, \Omega_X^{[1]}(pK_X))=0$.
By Serre duality of Cohen-Macaulay sheaves (\cite{KM98}*{Theorem 5.71}), we have $H^1(X, \Omega_X^{[1]}(-pK_X))\cong H^1(X, \Omega_X^{[1]}(pK_X))$.
Since $X$ has only hypersurface singularities, $\Omega^1_X$ is torsion-free by \cite{Lipman}*{Section 8 (1)}, and the natural map $\Omega_X^1\to \Omega_X^{[1]}$ is injective.
Since $\sO_X(-pK_X)$ is Cartier, we have an exact sequence
\[
0 \to \Omega^{1}_X\otimes \sO_X(-pK_X) \to \Omega^{[1]}_X(-pK_X) \to \mathcal{C} \to 0
\]
for some coherent sheaf $\mathcal{C}$ satisfying $\dim \Supp(\mathcal{C})=0$. 
Since $H^1(X,\mathcal{C})=0$, it suffices to show that 
\[
H^{1}(X, \Omega^{1}_X\otimes \sO_X(-pK_X))=0.
\]
In what follows, we divide the proof into the cases according to (1)--(4) in the proposition.

\textbf{The case (1):}
In this case, $X$ is a hypersurface of $P\coloneqq \mathbb{P}(1,1,2,3)$ of degree $6$ (\cite{BT22}*{Theorem 2.15}).
{By \cite{Mori75}*{Theorem 1.7}, the non-Gorenstein locus of $P$ is $\{[0:0:0:1], [0:0:1:0]\}$, and this locus coincides with the singular locus of $P$ (Remark \ref{r-toric}).
Thus, $X$ is contained in the smooth locus of $P$ since it is Gorenstein. 
We define invertible sheaves $\sO_X(n)$ by $\sO_P(n)|_X$ for all $n\in\Z$.}

By adjunction, we have $\omega_X=\sO_X(-1)$, and thus we aim to show that \[H^{1}(X, \Omega^1_X\otimes \sO_X(p))=0.\] 
By the conormal exact sequence, we have an exact sequence
\[
\sO_X(-X+p)=\sO_X(p-6)\to \Omega^1_{P}|_X\otimes \sO_X(p) \to \Omega^1_X\otimes \sO_X(p)\to 0.
\]
Since $\sO_X(p-6)$ is torsion-free and the first map is injective outside the singular points of $X$, we obtain an exact sequence
\[
0\to \sO_X(p-6)\to (\Omega^1_{P}\otimes \sO_P(p))|_X \to \Omega^1_X\otimes \sO_X(p)\to 0.
\]
Since $p\geq 7$, we have $H^2(X, \sO_X(p-6))=0$ by Serre duality, 
and hence it suffices to show that $H^1(X, (\Omega^1_{P}\otimes \sO_P(p))|_X)=0$.
We have an exact sequence
\[
\Omega^{[1]}_P(p-6)=\Omega^{[1]}_P(p)\otimes \sO_P(-6) \to \Omega^{[1]}_P(p) \to \Omega^{[1]}_P(p)|_X\to 0,
\]
Here, we obtain the first equality as follows:
\[
\Omega^{[1]}_P(p)\otimes \sO_P(-6)=(\Omega^{[1]}_P\otimes \sO_P(p))^{**}\otimes \sO_P(-6)=(\Omega^{[1]}_P\otimes \sO_P(p-6))^{**}=\Omega^{[1]}_P(p-6)
\]
since $\sO_P(-6)$ is Cartier.
In particular, the first term of the above exact sequence is torsion-free, and thus the first map is injective since it is injective outside the singular points of $P$.

Moreover, since $X$ is contained in the smooth locus of $P$, it follows that $\Omega^{[1]}_P(p)|_X=(\Omega^1_{P}\otimes \sO_P(p))|_X$.
Thus, we obtain an exact sequence
\[
0\to \Omega^{[1]}_P(p-6) \to \Omega^{[1]}_P(p) \to (\Omega^1_{P}\otimes \sO_P(p))|_X\to 0.
\]
By Bott vanishing on $P$ \cite{Fuj07}*{Corollary 1.3}, we have \[H^1(P, \Omega^{[1]}_P(p))=H^2(P, \Omega^{[1]}_P(p-6))=0\]
since $p\geq 7$.
Therefore, we obtain $H^1(X, (\Omega^1_{P}\otimes \sO_P(p))|_X)=0$.

\textbf{The case (2):}
In this case, $X$ is a hypersurface of $P\coloneqq \mathbb{P}(1,1,1,2)$ of degree $4$ (\cite{BT22}*{Theorem 2.15}).
{By \cite{Mori75}*{Theorem 1.7}, the non-Gorenstein locus of $P$ is $\{[0:0:0:1]\}$, and this locus coincides with the singular locus of $P$ (Remark \ref{r-toric}).
Thus, $X$ is contained in the smooth locus of $P$ since it is Gorenstein.}

By adjunction, we have $\omega_X=\sO_X(-1)$, and thus we aim to show that \[H^{1}(X, \Omega^1_X\otimes \sO_X(p))=0.\] 
As in the case (1), by the conormal exact sequence and the torsion-freeness of $\sO_X(p-4)$, we have an exact sequence
\[
    0\to \sO_X(p-4)\to \Omega^1_{P}|_X\otimes \sO_X(p) \to \Omega^1_X\otimes \sO_X(p)\to 0. 
\]
Since $p\geq 5$, we have $H^2(X, \sO_X(p-4))=0$, and it suffices to show that \[H^1(X, \Omega^1_{P}|_X\otimes \sO_X(p))=0.\]
As in the case (1), we have an exact sequence
\[
0\to \Omega^{[1]}_P(p-4)\to \Omega^{[1]}_P(p) \to \Omega^1_{P}|_X\otimes \sO_X(p)\to 0.
\]
By Bott vanishing \cite{Fuj07}*{Corollary 1.3}, we have \[H^1(P, \Omega^{[1]}_P(p))=H^2(P, \Omega^{[1]}_P(p-4))=0\]
since $p\geq 5$.
Therefore, we obtain $H^1(X, \Omega^1_{P}(p)|_X)=0$.

\textbf{The case (3):}
In this case, $X$ is a hypersurface of $P\coloneqq \mathbb{P}^3$ of degree $3$ (\cite{BT22}*{Theorem 2.15}).
By adjunction, we have $\omega_X=\sO_X(-1)$, and thus we aim to show that $H^{1}(X, \Omega^1_X(p))=0$. 
By the conormal exact sequence and the torsion-freeness of $\sO_X(p-3)$, we have an exact sequence
\[
0\to \sO_X(p-3)\to \Omega^1_{P}(p)|_X \to \Omega^1_X(p)\to 0.
\]
Since $p\geq 3$, we have $H^2(X, \sO_X(p-3))=0$, and it suffices to show that $H^1(X, \Omega^1_{P}(p)|_X)=0$.
We have an exact sequence
\[
0 \to \Omega^{1}_P(p-3) \to \Omega^{1}_P(p) \to \Omega^{1}_P(p)|_X\to 0.
\]
By Bott vanishing \cite{Fuj07}*{Corollary 1.3}, we have $H^1(P, \Omega^{1}_P(p))=0$.
By \cite{Totaro(Fano)}*{Proposition 1.3}, we also have $H^2(P, \Omega^{1}_P(p-3))=0$ since $p\geq 3$.
Therefore, we obtain $H^1(X, \Omega^1_{P}(p)|_X)=0$.

\textbf{The case (4):}
In this case, $X$ is a complete intersection of two 
quadric hypersurfaces $Q$ and $Q'$ of $P\coloneqq \mathbb{P}^4$ (\cite{BT22}*{Theorem 2.15}).
By adjunction, we have $\omega_X=\sO_X(-1)$, and thus we aim to show that $H^{1}(X, \Omega^1_X(p))=0$. 
By the conormal exact sequence and the torsion-freeness of $\sO_X(p-2)$, we have an exact sequence
\[
0\to \sO_X(p-2)\to (\Omega^1_{Q}\otimes \sO_{Q}(p))|_X \to \Omega^1_X\otimes \sO_X(p)\to 0.
\]
Since $p\geq 2$, we have $H^2(X, \sO_X(p-2))=0$, and hence it suffices to show that $H^1(X, (\Omega^1_{Q}\otimes \sO_{Q}(p))|_X)=0$.

We define invertible sheaves $\sO_Q(n)$ by $\sO_{P}(n)\otimes \sO_{Q}$ for all $n\in\Z$.
We have an exact sequence
\[
\Omega^{1}_Q\otimes \sO_Q(p-2) \to \Omega^{1}_Q\otimes \sO_Q(p) \to (\Omega^1_{Q}\otimes \sO_{Q}(p))|_X \to 0.
\]
Since $X$ is regular in codimension one, $\Omega^1_Q$ is torsion-free by \cite{Lipman}*{Section 8 (1)}.
Thus, we have an exact sequence
\[
0\to \Omega^{1}_Q\otimes \sO_Q(p-2) \to \Omega^{1}_Q\otimes \sO_Q(p) \to (\Omega^1_{Q}\otimes \sO_{Q}(p))|_X \to 0.
\]
Therefore, it suffices to show that
\[H^1(Q, \Omega^1_Q\otimes \sO_Q(p))=0\,\,\,\text{and}\,\,\,H^2(Q, \Omega^1_Q\otimes \sO_Q(p-2))=0,\]
and in particular, the following claim finishes the proof of the case (4):

\begin{claim*}
    We have
\begin{enumerate}
    \item[(i)] $H^1(Q, \Omega^1_{Q}\otimes \sO_Q(n))=0$ for every $n\in \Z \setminus \{0\}$ and
    \item[(ii)] $H^2(Q, \Omega^1_{Q}\otimes \sO_Q(n))=0$ for every $n\in\Z$.
\end{enumerate}
\end{claim*}
We have
\begin{enumerate}
\item[(a)] $H^i(P, \Omega^1_{P}(n))=0$ for every $n \in \Z \setminus \{0\}$ and $i \in \{1, 2, 3\}$, and
\item[(b)] $H^j(P, \Omega^1_P(n))=0$ for every $n \in \Z$ and $j \in \{2, 3\}$.  
\end{enumerate}
Indeed, (a) follows from Bott vanishing and Serre duality. 
Then (a), together with \cite{Totaro(Fano)}*{Proposition 1.3}, implies (b). 
By the following exact sequence
\[
0\to \Omega^1_{P}(n-2)\to \Omega^1_{P}(n) \to \Omega^1_{P}(n)|_{Q} \to 0, 
\]
we get 
\begin{enumerate}
    \item[(i)'] $H^1(Q, \Omega^1_{P}(n)|_{Q})=0$ for every $n\in \Z\setminus \{0\}$, and 
    \item[(ii)'] $H^2(Q, \Omega^1_{P}(n)|_{Q})=0$ for every $n\in\Z$.
\end{enumerate}
By the conormal exact sequence and the torsion-freeness of $\sO_{P}(n-2)$, we have an exact sequence
\[
0\to \sO_{P}(n-2)\to \Omega^1_{P}(n)|_{Q} \to \Omega^1_Q\otimes \sO_Q(n) \to 0
\]
for every $n\in\ Z$. 
Since $H^2(P, \sO_P(n))=H^3(P, \sO_P(n))=0$ for every $n\in\Z$, we have the claim.
\end{proof}

\begin{rem}\label{r-toric}
Take positive integers $q_1, q_2, q_3$ such that ${\rm gcd}(q_1, q_2, q_3)=1$. 
Set $P\coloneqq\P(1, q_1, q_2, q_3)$. 
Then it is well known (cf.\ \cite{Ful93}*{Section 2.2}) that $P$ coincides with the projective $\Q$-factorial toric threefold associated to the fan in $\R^3$ that is generated by four rays $\R u, \R e_1, \R e_2, \R e_3$, where $e_1, e_2, e_3$ is the standard basis of $\Z^3$ and 
\[
u \coloneqq -(q_1e_1+q_2e_2+q_3e_3). 
\]
In the above proof, we have used the results (1) and (2). 
\begin{enumerate}
\item 
$\P(1, 1, 2, 3)$ has exactly two singular points, 
which corresponds to the cones $\R u + \R e_1 + \R e_2$ and  $\R u + \R e_1 + \R e_3$ \cite{CLS11}*{Theorem 1.3.12}. 
\item $\P(1, 1, 1, 2)$ has a unique singular point, which corresponds to the cone $\R u + \R e_1 + \R e_2$ \cite{CLS11}*{Theorem 1.3.12}. 
\end{enumerate}
\end{rem}

From now on, we focus on the case where $K_X^2 \geq 5$.

\begin{prop}\label{p-distinct}
The following assertions hold.\,
\begin{enumerate}
\item 
Fix an integer $m$ satisfying $1 \leq m \leq 5$. 
Let $P_1,\ldots, P_m$ be distinct points on $\P^2$ such that
 the blowup $X$ of $\P^2$ along  $\{P_1, \ldots, P_m\}$ is a weak del Pezzo surface. 
Then $X$ is $F$-split.  
\item 
Fix an integer $n$ satisfying $1 \leq n \leq 4$. 
Let $Q_1,\ldots, Q_n$ be distinct points on $\P^1 \times \P^1$ such that
 the blowup $X$ of $\P^1 \times \P^1$ along  $\{Q_1,\ldots, Q_n\}$ is a weak del Pezzo surface. 
Then $X$ is $F$-split.  
\end{enumerate}
\end{prop}

\begin{proof}
Let us show (1). 
In what follows, we only treat the case when $m=5$, as otherwise the problem is easier. 
Let $L$ (resp.~$L'$) be the line on $\P^2$ passing through $P_1$ and $P_2$ (resp.~$P_3$ and $P_4$). 
Since $X$ is weak del Pezzo, we obtain $L \neq L'$. 
Pick a general line $L''$ on $\P^2$ passing through $P_5$. 
Then $L+L'+L''$ is simple normal crossing. 
Therefore, $(\P^2, L+L'+L'')$ is $F$-split (Example \ref{e-P2-P1P1}(2)), which implies that so is $X$  (Proposition \ref{p-F-split-blowup}). 
Thus (1) holds. 
The proof of (2) is similar to that of (1). Indeed, 
for each projection $\pr_i \colon \P^1 \times \P^1 \to \P^1$, 
it is enough to take two fibers $F_i$ and $F'_i$ such that 
$F_1 \cup F'_1 \cup F_2 \cup F'_2$ contains $\{Q_1,\ldots, Q_n\}$  (Example \ref{e-P2-P1P1}(3)). 
\end{proof}

\begin{prop}\label{prop:degree 5}
Let $X$ be a smooth weak del Pezzo surface satisfying $K_X^2 \geq 5$. 
Then $X$ is $F$-split. 
\end{prop}

\begin{proof}
By \cite{Dol12}*{Theorem 8.1.15}, 
we may assume that there is a birational morphism $f\colon X \to \P^2$. 
In what follows, we only treat the case when $K_X^2 = 5$, 
as the other cases are simpler. 
There are the following five cases \cite{Dol12}*{Section 8.5}. 
\begin{enumerate}
\renewcommand{\labelenumi}{(\roman{enumi})}
\item $P, Q, R, S$. 
\item $P' \succ P, Q, R$. 
\item $P' \succ P, Q' \succ Q$. 
\item $P'' \succ P' \succ P, Q$.
\item $P'' \succ P' \succ P  \succ Q$. 
\end{enumerate}
For the definition of $P' \succ P$, 
we refer to  \cite{Dol12}*{Section 7.3.2}. 
For example, in the case (iii), 
we have $X =Y'' \to Y' \to Y = \P^2$, where 
$Y' = \Bl_{P \amalg Q} Y$, $Y'' = \Bl_{P' \amalg Q'}$, 
and $P'$ and $Q'$ are points on $Y'$ lying over $P$ and $Q$, respectively.

The case (i) has been settled in Proposition \ref{p-distinct}. 
As for (ii), take the line $L_1 \coloneqq \overline{PQ}$ passing through $P$ and $Q$. 
Let $L_2$ and $L_3$ be general lines passing through $P$ and $R$, respectively. 
Then $(\P^2, L_1+L_2+L_3)$ is $F$-split  (Example \ref{e-P2-P1P1}(2)). 
Since $\Delta$ is effective for the divisor $\Delta$ defined by 
$K_X +\Delta =f^*(K_{\P^2}+L_1+L_2+L_3)$, it follows that
$X$ is $F$-split (Proposition \ref{p-F-split-blowup}). 
Similarly, 
(iii) is settled by taking the line $L_1 \coloneqq \overline{PQ}$ and general lines $L_2$ and $L_3$ passing through $P$ and $Q$, respectively. 

Let us treat the case (iv). 
In this case, we have $X = Y''' \to Y'' \to Y' \to Y=\P^2$, 
where $Y' \coloneqq \Bl_{P\amalg Q} Y, Y'' \coloneqq \Bl_{P'} Y'$, 
and $Y''' \coloneqq \Bl_{P''} Y''$. 
Let $L_1$ be the line on $Y=\P^2$ 
such that $P \in L_1$ and $P' \in L'_1$ 
for the proper transform $L'_1$ of $L_1$ on $Y'$. 
Let $L_2$ and $L_3$ be general lines passing through $P$ and $Q$, respectively. 
Then we can check that the divisor $\Delta$ defined by $K_X +\Delta =f^*(K_{\P^2}+L_1 +L_2 +L_3)$ is 
effective. 
Since $(\P^2, L_1 +L_2 +L_3)$ is $F$-split  (Example \ref{e-P2-P1P1}(2)), so is $X$. 
This completes the proof for the case (iv). 

\medskip

Let us consider the case (v). 
In this case, we apply a similar method to that of (iv) after replacing $\P^2$ by $\F_1$. 
We have a sequence of one-point blowups: 
\[
f\colon X = Y''' \to Y'' \to Y' \to Y=\F_1 \to Z =\P^2,
\]
{where $Y' \coloneqq \Bl_{P} Y, Y'' \coloneqq \Bl_{P'} Y'$, and $Y''' \coloneqq \Bl_{P''} Y''$.}
For the $(-1)$-curve $C$ on $Y$, we have $P\in C$. 
It is well known that there is another section $\widetilde C$ of the $\P^1$-bundle $\pi \colon Y=\F_1 \to B=\P^1$ such that $C \cap \widetilde{C} = \emptyset$ and $\widetilde{C}^2 =1$. 
{Let $F$ is a fiber of $\pi$. Since $(K_Y+C+\widetilde{C})\cdot F=0$, there exists $n\in \Z$ such that $K_Y+C+\widetilde{C}\sim nF$. 
Then $n=C\cdot nF =C\cdot (K_Y+C+\widetilde{C})=-2$.} 
Since the proper transform $C'$ of $C$ on $Y'$ satisfies $C'^2=-2$, 
we obtain 
\begin{enumerate}
    \item[(1)] $P' \not\in C'$, 
\end{enumerate} 
as otherwise, the proper transform $C''$ of $C'$ on $Y''$ would satisfy 
$C''^2 =-3$, which contradicts the fact that $Y''$ is weak del Pezzo. 

We now treat the case when $P' \in F'_P$, 
where $F'_P$  denotes the proper transform 
of the fiber $F_P$ of $\pi \colon Y =\F_1 \to \P^1$ passing through $P$.
Let $\widetilde F$ be a general fiber of $\pi$. As we have seen above, $K_Y+{\widetilde F}+ C+F_P +\widetilde C \sim 0$.
Since $\widetilde F$ is nef and $\F_1$ is toric, we obtain $H^1(\F_1, \sO_{F_1}(\widetilde F))=0$ by \cite{Totaro(Fano)}*{Proposition 1.3}. Moreover, $(\widetilde F, (C+F_P +\widetilde C)|_{\widetilde F})=( \widetilde{F}, C|_{\widetilde F}+\widetilde C|_{\widetilde F})$ is $F$-split (Example \ref{e-P2-P1P1} (1)).
Thus, $(Y, C+F_P +\widetilde C + \widetilde F)$ is $F$-split (Proposition \ref{p-F-split-IOA}), 
which implies that $X$ is $F$-split (Proposition \ref{p-F-split-blowup}).
In what follows, we assume that 
\begin{enumerate}
\item[(2)] $P' \not\in F'_P$ for the proper transform $F'_P$ 
of the fiber $F_P$ of $\pi : Y =\F_1 \to \P^1$ passing through $P$. 
\end{enumerate}

\begin{claim*}
There is a section $D$ of $\pi \colon Y=\F_1 \to \P^1$ such that 
\begin{enumerate}
\item[(a)] $D \sim \widetilde C+F$, 
\item[(b)] $P \in D$, and 
\item[(c)] $P' \in D'$ for the proper transform $D'$ of $D$ on $Y'$. 
\end{enumerate}
\end{claim*}

\begin{proof}[Proof of Claim]
Since $\widetilde C+F$ is an ample Cartier divisor on $Y=\F_1$, 
it follows that $|\widetilde C+F|$ is very ample \cite{Har77}*{Ch.~V, Corollary 2.18}. 
Then there is an effective Cartier divisor $D$ on $Y = \F_1$ satisfying (a)--(d).  
\begin{enumerate}
\item[(d)] $D$ is smooth at $P$. 
\end{enumerate}
In fact, since $|\widetilde C+F|$ is very ample, the elements of $H^0(Y, \sO_Y(\widetilde C+F))$ separate tangent vectors.
Let $s_{P'}\in \m_P/\m^2_{P}$ is an element that corresponds to $P'$.
We take $D$ as a divisor of zeros of a global section $s\in H^0(Y, \sO_Y(\widetilde C+F))$ that maps to $s_{P'}\in \m_P/\m^2_{P}$. 
Then (a)--(c) are satisfied. Since $s \in \m_P/\m^2_{P}$ is non-zero, the divisor $D$ is smooth at $P$, i.e., (d) is satisfied.
Since $D \cdot F = (\widetilde C + F) \cdot F =1$, 
we can write 
$D= D_0 + F_1 + \cdots +F_r$, 
where 
$r \geq 0$, 
$D_0$ is a section of $\pi \colon Y = \F_1 \to \P^1$, and 
each $F_i$ is a fiber of $\pi$. 
It suffices to prove $r=0$. 
Suppose $r>0$. 
The following holds: 
\begin{equation}\label{e1-quintic-dP}
D_0 \cdot C +r = (D_0+F_1 + \cdots +F_r) \cdot C =D \cdot C = (\widetilde C +F) \cdot  C =1. 
\end{equation}

We now treat the case when $D_0 \neq C$. 
In this case, $D_0 \cdot C \geq 0$ and (\ref{e1-quintic-dP}) imply $r=1$ and $D_0 \cdot C =0$. 
Hence we get $D_0 \cap C = \emptyset$. 
Since $P\in C$, we have $P\notin D_0$.  By (b), we obtain $P\in D=D_0+F_1$, and thus $P\in F_1$.
Hence $D = D_0 +F_P$, 
where $F_P$ denotes the fiber passing through $P$. 
{Since $P'\notin C'$, we obtain $P'\in D'\setminus C'\subset F'_P$.}
This contradicts (2).  

Hence we may assume that $D_0 =C$. 
We then get $D = C +F_1 + \cdots +F_r$. 
Since $P\in C$, we obtain $P \not\in F_1 \cup \cdots \cup F_r$ by (d). 
{Then $P' \not\in F'_1 \cup \cdots \cup F'_r$, where $F'_1, F'_2,\ldots, F'_r$ are proper transforms of $F_1, F_2,\ldots, F_r$ on $Y'$, and thus we obtain $P'\in D'\setminus \{F'_1 \cup \cdots \cup F'_r\}\subset C'$ by (c).}
This contradicts (1). 
This completes the proof of Claim. 
\end{proof}

We have $C \cdot D = 1$. 
Hence $C \cap D = P$ and $C+ D$ is a simple normal crossing divisor. 
Since both $C$ and $D$ are sections of $\pi \colon Y = \F_1 \to \P^1$, it follows that $C+D+\widetilde F$ is still simple normal crossing for a general fiber $\widetilde F$ of $\pi$. 
Then we see that $(Y, C+D+\widetilde F)$ is $F$-split (Proposition \ref{p-F-split-IOA}, Example \ref{e-P2-P1P1}(1)). 
Since  
the divisor $\Delta$ defined by $K_X +\Delta  = f^*(K_Y+C+D+\widetilde F)$ is effective, $X$ is $F$-split  (Proposition \ref{p-F-split-blowup}). 
This completes the proof of Proposition \ref{prop:degree 5}. 
\end{proof}

\begin{proof}[Proof of Theorem \ref{Intro:F-split}]
If $K_X^2 \leq 4$ (resp.\ $K_X^2 \geq 5$), then 
the assertion follows from Proposition \ref{prop:degree 1--4} 
(resp.\ Proposition \ref{prop:degree 5}).
\end{proof}

\subsection{Global \texorpdfstring{$F$}--regularity: Proof of Theorem \ref{Intro:main}}

In this subsection, we deduce Theorem \ref{Intro:main} from Theorem \ref{Intro:F-split}.

\begin{lem}\label{lem:F-reg-sing}
    Let $X$ be a canonical del Pezzo surface.
    Suppose that one of the following holds.
\begin{enumerate}
    \item $p>5$
    \item $K_X^2\geq 2$ and $p>3$.
    \item $K_X^2\geq 4$ and $p>2$.
    \item $K_X^2\geq 5$.
\end{enumerate}
Then $X$ is strongly $F$-regular. 
\end{lem}
\begin{rem}
    Combining \cite{KN2}*{Table 1} and \cite{Kawakami-Takamatsu}*{Table 1}, we can see that the assumption of $p$ is optimal for each degree.
\end{rem}

\begin{proof}
(1) follows from \cite{Hara(two-dim)}*{Theorem 1.1}.
In what follows, let $f\colon Y\to X$ be the minimal resolution.
    
We prove (2).
Since $K_Y^2=K_X^2\geq 2$, we have $\rho(Y) = 10 -K_Y^2 \leq 8$. Thus, the number of the $(-2)$-curves contracted by $f$ is at most $8-\rho(X)\leq 7$.
Therefore, $X$ does not have canonical singularities of $E_8$-type. 
Then $X$ is strongly $F$-regular by \cite{Hara(two-dim)}*{Theorem 1.1} since $p>3$ (see also \cite{Kawakami-Takamatsu}*{Table 1}).

    Next, we prove (3).
    Since $K_Y^2=K_X^2\geq 4$, we have $\rho(Y) = 10 -K_Y^2\leq 6$.
    Thus, the number of the $(-2)$-curves contracted by $f$ is at most $6-\rho(X)\leq 5$. 
    Therefore, $X$ does not have canonical singularities of $E$-type.
    Then $X$ is strongly $F$-regular by \cite{Hara(two-dim)}*{Theorem 1.1} since $p>2$ (see also \cite{Kawakami-Takamatsu}*{Table 1}).

    Finally, we prove (4).
    Since $K_Y^2=K_X^2\geq 5$, we have $\rho(Y) = 10 -K_Y^2\leq 5$.
    Thus, the number of the $(-2)$-curves contracted by $f$ is at most $5-\rho(X)\leq 4$.
    If $\rho(X)\geq 2$, then $X$ has only $A$-type singularities, which are strongly $F$-regular \cite{Hara(two-dim)}*{Theorem 1.1}.
    If $\rho(X)=1$, then $X$ has only $A$-type singularities by \cite{KN2}*{Theorem 1.1}.
\end{proof}

\begin{lem}\label{lem:Kawakami-Totaro}
    Let $X$ be a normal projective variety such that $-K_X$ is an ample $\Q$-Cartier $\Z$-divisor.
    Suppose that $X$ is strongly $F$-regular. 
    If $X$ is $F$-split, then $X$ is globally $F$-regular.
\end{lem}
\begin{proof}
    See the proof of \cite{Kawakami-Totaro}*{Lemma 6.4}. 
\end{proof}

\begin{proof}[Proof of Theorem \ref{Intro:main}]
Let $X$ be as in the statement of  Theorem \ref{Intro:main}. 
Taking the anti-canonical model of $X$, we may assume that $-K_X$ is ample (Corollary \ref{c-F-split-blowup}). 
Then, by Lemma \ref{lem:Kawakami-Totaro}, it is enough to prove that $X$ is strongly $F$-regular and $F$-split, which follow from Lemmas \ref{lem:F-reg-sing} and Theorem \ref{Intro:F-split}, 
respectively. 
\end{proof}


\end{document}